\documentclass[12pt]{amsart}
\usepackage[latin1]{inputenc}
\usepackage[T1]{fontenc}
\usepackage[french]{babel}
\usepackage{amssymb,amsfonts,latexsym}
\usepackage[arrow,matrix,curve]{xy}
\usepackage[mathscr]{eucal}
\usepackage{a4wide}
\usepackage{enumerate}

\font\trmbi=cmtcsc10 scaled 1200
\font\srmbi=cmtcsc10 scaled 840         
\font\ssrmbi=cmtcsc10 scaled 600
\font\teuf=eufm10 scaled 1200
\font\seuf=eufm7 scaled 1200            
\font\sseuf=eufm5 scaled 1200
\font\tenscr=rsfs10  scaled 1200
\font\sevenscr=rsfs7 scaled 1200        
\font\fivescr=rsfs5  scaled 1200

\skewchar\tenscr='177 \skewchar\sevenscr='177 \skewchar\fivescr='177

\newfam\majfam
        \textfont\majfam=\trmbi
        \scriptfont\majfam=\srmbi
        \scriptscriptfont\majfam=\ssrmbi

\newfam\euffam
        \textfont\euffam=\teuf
        \scriptfont\euffam=\seuf
        \scriptscriptfont\euffam=\sseuf

\newfam\scrfam 
        \textfont\scrfam=\tenscr 
        \scriptfont\scrfam=\sevenscr
        \scriptscriptfont\scrfam=\fivescr

\def \goth{\fam\euffam}

\def \scr{\fam\scrfam}

\def \SSS{\scriptscriptstyle}

\def \iso{\hbox{$\,\vbox{\hbox to
18pt{\hfill$\sim$\hfill}\kern-11pt\hbox{$\longrightarrow$}}\,$}}
\def \liso{\hbox{$\,\vbox{\hbox to
18pt{\hfill$\sim$\hfill}\kern-11pt\hbox{$\longleftarrow$}}\,$}}
\def \doublefleche{\vbox
to5.5025pt{\hbox{$\longrightarrow$}\kern-10.75pt\hbox{$\longrightarrow$}}}
\def \subsetneq{\;\lower1.66667pt\vbox
to9.60166pt{\vskip0.3pt\hbox{$\subseteq$}\kern-12.83pt\hbox{$\;\SSS/\;$}}\;}

\def \subsetneq{\;\lower1.66667pt\vbox
to9.60166pt{\vskip0.3pt\hbox{$\subseteq$}\kern-12.83pt\hbox{$\;\SSS/\;$}}\;}
\def \iso{\buildrel\sim\over\longrightarrow}

\def \HM{{\mathbb{H}}}

\def \NM{{\mathbb{N}}}

\def \QM{{\mathbb{Q}}}

\def \ACC{{\scr A}}

\def \CCC{{\scr C}}

\def \ECC{{\scr E}}
\def \FCC{{\scr F}}

\def \HCC{{\scr H}}
\def \ICC{{\scr I}}
\def \JCC{{\scr J}}

\def \LCC{{\scr L}}
\def \MCC{{\scr M}}
\def \NCC{{\scr N}}
\def \OCC{{\scr O}}

\def \VCC{{\scr V}}

\def \XCC{{\scr X}}
\def \YCC{{\scr Y}}
\def \ZCC{{\scr Z}}

\def \UG{{\goth U}}

\def \leq{\leqslant}
\def \geq{\geqslant}

\newcommand{\rqe}{\noindent {\sc Remarque : }\rm}
\newcommand{\notation}{\noindent {\sc Notations : }\rm}

\newcommand {\dem}{\noindent {\sc D{\'e}monstration : } \rm }
\newcommand {\findem}{\hfill$\Box$\par\medskip}
\newcount\exno
\newcommand{\exo}[1]{\advance\exno by1 \medskip
                \hangindent=\parindent\hangafter=1
                \noindent{\trmbi Exercice}
                \hbox{ \arabic{chapter}.\the\exno.}\hfill \\ 
               \rm #1\par\medskip}

\renewcommand{\thesubsection}{\arabic{section}.\arabic{subsection}}

\newtheorem{prop}{Proposition}[subsection]
\newtheorem{lemme}[prop]{Lemme}
\newtheorem{Defi}[prop]{D{\'e}finition}
\def\moncompteur{\thesubsection.\arabic{prop}}
\def \texte #1 {\addtocounter{prop}1\noindent(\moncompteur) {\it #1}\\}
\newtheorem{theo}[prop]{Th{\'e}or{\`e}me}
\newtheorem{coro}[prop]{Corollaire}

\newtheorem{theonn}{Th{\'e}or{\`e}me}

\makeatletter

\@addtoreset{equation}{section}
\makeatother

\def \limind#1{\lim\limits_{\displaystyle\longrightarrow\atop {#1}}}

\def \dim{\mathop{\rm dim}\nolimits}
\def \spec#1{\hbox{\rm Spec}\,(#1)}
\def \spf#1{\hbox{\rm Spf}\,(#1)}
\def \spm#1{\hbox{\rm Spm}\,(#1)}

\def \hrc#1{H^{#1}_{c,rig}}
\def \hrz#1#2{H^{#1}_{#2,rig}}

\def \ov#1{\overline{#1}}
\def \surfleche#1{\buildrel#1\over\longrightarrow}
\def \dlog#1{\frac{d#1}{#1}}

\title{Classes de cycles en cohomologie rigide}
\author{Denis Petrequin }

\begin{document}
\maketitle
\tableofcontents
\bibliographystyle{plain}

\section*{Introduction}

\def \cl{{\rm cl}}

Nous allons montrer qu'il existe une th{\'e}orie
des classes fondamentales et des classes de cycles en cohomologie rigide. On pourrait s'appuyer sur
\cite{SGA6} pour d{\'e}duire directement pour les
vari{\'e}t{\'e}s quasi-projectives l'existence des classes de cycles {\`a} partir des
classes de Chern \cite{Pet1}. Cependant, nous allons donner une construction
directe qui a le m{\'e}rite, en plus de pouvoir traiter le cas des vari{\'e}t{\'e}s non n{\'e}cessairement quasi-projectives, d'{\^e}tre plus explicite. Dans le cas des vari\'et\'es singuli\`eres, les classes de cycles se construisent naturellement dans une {\it homologie rigide} que nous d\'efinirons. Nous v\'erifierons alors le th\'eor\`eme suivant :

\begin{theonn}
Soient $k$ un corps, $C$ un anneau de Cohen pour $k$, $K$ le corps des fractions de $C$ et $\VCC$ la cat\'egorie des $k$-vari\'et\'es. Le couple cohomologie rigide - homologie rigide forme une th\'eorie de dualit\'e de Poincar\'e au sens de Bloch-Ogus \cite{Jann, Blo-Ogu}. Plus pr\'ecis\'ement on a :
\begin{enumerate}
\item \label{morph} $H^i_{Y,rig}(X/K)$ est contravariant par rapport aux carr\'es cart\'esiens 
$$\xymatrix{ Y \ar@{^{(}->}[r] \ar[d] & X \ar[d] \\ Y' \ar@{^{(}->}[r] & X'.}$$
\item $H_i^{rig}(X/K)$ est contravariant par rapport aux immersions ouvertes et covariant par rapport aux morphismes propres.
\item \label{sel} Pour $X$ une $k$-vari\'et\'e et $Z \subset Y \subset X$ des sous-sch\'emas ferm\'es, il existe une suite exacte longue 
$$\cdots \to H^i_{Z,rig}(X/K) \to  H^i_{Y,rig}(X/K) \to H^i_{Y-Z,rig}(X-Z/K) \to H^{i+1}_{Z,rig}(X/K) \to \cdots,$$
fonctorielle par rapport aux morphismes de fonctorialit\'e de \ref{morph}.
\item {\bf (excision)}  Pour tout $X$, $U$ un ouvert de $X$ et $Z$ un ferm\'e de $X$ tel que $Z \subset U$, le morphisme $H^i_{Z, rig}(X/K) \to H^i_{Z,rig}(U/K)$ est un isomorphisme.
\item \label{1636} Si on a le diagramme cart\'esien suivant
$$\xymatrix{ X' \ar@{^{(}->}[r]^{\beta} \ar_{g}[d] & X \ar^{f}[d] \\ Y' \ar@{^{(}->}[r]_{\alpha} & Y}$$
avec $\alpha$ et $\beta$ des immersions ouvertes et $f$ et $g$ des morphismes propres, le diagramme suivant est commutatif
$$\xymatrix{H_i^{rig}(X/K) \ar[r]^{\beta^*} \ar_{f_*}[d] & H_i^{rig}(X'/K) \ar^{g_*}[d] \\ H_i^{rig}(Y/K) \ar[r]_{\alpha^*} &  H_i^{rig}(Y'/K).}$$
\item Si $i:Y \hookrightarrow X$ est une immersion ferm\'ee et si $\alpha : X-Y \hookrightarrow X$ est l'immersion ouverte compl\'ementaire, il existe une suite exacte longue
$$\cdots \to H_{i+1}^{rig}(X-Y/K) \to H_{i}^{rig}(Y/K) \surfleche{i_*} H_{i}^{rig}(X/K) \surfleche{\alpha^*} H_{i}^{rig}(X-Y/K) \to H_{i-1}^{rig}(Y/K) \to \cdots $$
qui est fonctorielle par rapport aux morphismes propres.
\item {\bf (cap-produit)} Il existe un cap-produit pour tout $X$ et tout ferm\'e $Y \hookrightarrow X$ :
$$H_i^{rig}(X/K) \otimes H^j_{Y,rig}(X/K) \surfleche{\cap} H_{i-j}^{rig}(Y/K).$$
\item {\bf (formule de projection)} Pour tout diagramme cart\'esien 
$$\xymatrix{ Y' \ar@{^{(}->}[r]^{i'} \ar[d]_{f'} & X' \ar[d]^{f} \\ Y \ar@{^{(}->}[r]_{i} & X,}$$
o\`u $i$ et $i'$ sont des immersions ferm\'ees et $f$ et $f'$ des morphismes propres, et tous \'el\'ements $x' \in H_i^{rig}(X'/K)$ et $x \in H^j_{Y,rig}(X/K)$ on a :
$$f_*(x')\cap x = f'_*(x' \cap f^*(x)).$$

\item {\bf (classe fondamentale)} Pour tout vari\'et\'e irr\'eductible $X$ de dimension $d$, il existe une classe fondamentale
$$\eta_X \in H_{2d}^{rig}(X/K),$$
qui est fonctorielle par rapport aux immersions ouvertes.

\item {\bf (dualit\'e de Poincar\'e)} Si $X$ est une $k$-vari\'et\'e irr\'eductible et lisse de dimension $d$, et si $Y \hookrightarrow X$ est une immersion ferm\'ee, le morphisme
$$\xymatrix{ H_{Y,rig}^{2d-i}(X/K) \ar[rr]^{\eta_X \cap}&  & H_{i}^{rig}(Y/K)}$$ est un isomorphisme.
\item La dualit\'e de Poincar\'e est compatible aux suites exactes longues de \ref{sel}.
\end{enumerate}
\end{theonn} 

Par la suite, nous d\'efinirons les classes de cycles. Nous montrerons, \`a l'aide d'une comparaison avec la classe d'un diviseur qu'elles passent \`a l'\'equivalence rationnelle au sens \cite{Fult2}. 

On montrera alors que, en se limitant aux vari\'et\'es quasi-projectives, les classes de cycles d\'efinies pr\'ec\'edemment sont compatibles \`a la th\'eorie de l'intersection. 

On finira en \'enon\c{c}ant des cons\'equences de notre formalisme : th\'eor\`eme de Riemman-Roch, formule de self-intersection.

Cet article est la deuxi\`eme partie - avec quelques modifications - de ma th\`ese de doctorat \cite{moi}. Je tiens a remercier P.Berthelot qui m'a encadr\'e lors de ce travail.

\notation Dans tout cet article $k$ d\'esigne un corps de caract\'eristique $p>0$. On appelle $k$-vari\'et\'e un $\spec{k}$-sch\'ema s\'epar\'e de type fini. De plus, quand nous n'aurons pas \`a nous occuper du corps de base, nous noterons 

$$H^i_{Y,rig}(X) \hbox{ (resp. }H_i^{rig}(X) ) \hbox{ pour } H^i_{Y,rig}(X/K) \hbox{ (resp. }H_i^{rig}(X/K) ).$$

\section{D\'emonstration du th\'eor\`eme}
Nous allons d\'emontrer le th\'eor\`eme pr\'ec\'edent. Nous commencerons par rappeller certaines propri\'et\'es de la cohomologie rigide. Ensuite nous d\'efinirons l'homologie rigide et la classe fondamentale. Nous d\'emontrerons alors la formule de projection. Toutes les autres assertions d\'ecoulent directement des r\'esultats connus sur la cohomologie rigide \cite{Berth3, Berth7, Berth8}
\subsection{Rappels de cohomologie rigide}
Soient $X$ une $k$-vari{\'e}t{\'e} et
$Z$ un sous-sch{\'e}ma ferm{\'e} int{\`e}gre de $X$. On supposera que $X$ est {\'e}quidimensionel de dimension $n$ et que $Z$ est de codimension $r$. On notera alors $d = n-r$ la dimension de $Z$.

Rappelons pour commencer deux des principaux r{\'e}sultats de cohomologie
rigide. On trouvera leurs d{\'e}monstrations dans \cite{Berth7} et
\cite{Berth3} respectivement.

\begin{theo}[Dualit{\'e} de Poincar{\'e}]

 Avec les notations pr{\'e}c{\'e}dentes, il existe une application 
trace 
canonique $\hrc{2n}(X) \to K$ qui, compos{\'e}e avec la 
multiplication, induit des applications $\hrz{i}{Z}(X) \times
\hrc{2n-i}(Z) \to K$. De plus, si $X$ est lisse, c'est un accouplement 
parfait.
\end{theo}

\begin{theo}[Puret{\'e}]
Avec les notations pr{\'e}c{\'e}dentes, si on suppose que $X$ est lisse, on a pour tout $i<2r$ :
$$\hrz{i}{Z}(X) = 0.$$
De plus, la dimension sur $K$ de l'espace $\hrz{2d}{Z}(X)$ est \'egale au nombre de composantes irr\'eductibles g\'eom\'etriques de $Z$.
\end{theo}
\subsection {Homologie rigide}
\label{homol}
On se donne une $k$-vari{\'e}t{\'e} $X$, on pose
$$H_i^{rig}(X) := H^i_{c,rig}(X)^\vee.$$

Supposons qu'il existe une $k$-vari{\'e}t{\'e}  $M$ lisse sur $k$ et une immersion ferm{\'e}e $X \hookrightarrow M$.
En notant $N$ la dimension de $M$, on a gr{\^a}ce {\`a} la dualit{\'e} de Poincar{\'e} ($M$ est lisse) :
$$H_i^{rig}(X) \iso H^{2N-i}_{X,rig}(M).$$

\rqe si $X$ est lisse de dimension $n$, on a $H_i^{rig}(X) \iso H^{2n-i}_{rig}(X)$.

Regardons maintenant la variance de l'homologie rigide.
\begin{prop}
L'homologie rigide ainsi d{\'e}finie est contravariante par rapport aux
immersions ouvertes et covariante vis-{\`a}-vis des morphismes propres.
\end{prop}

\dem
Soit $f: X \to Y$ un morphisme propre. On sait que la cohomologie {\`a}
support compact est contravariante. Il existe donc pour tout $i$ :
$$f^* : H^i_{c,rig}(Y) \to H^i_{c,rig}(X).$$
En prenant la transpos{\'e}e on obtient le morphisme de fonctorialit\'e voulu :
$$f_* : H_{i}^{rig}(X) \to H_{i}^{rig}(Y).$$

De plus, on sait aussi que la cohomologie rigide {\`a} support compact
est covariante par rapport aux immersions ouvertes. Le m{\^e}me
raisonnement nous permet de construire le morphisme de fonctorialit{\'e} contravariante pour l'homologie rigide.

\findem

Nous allons \'enoncer quelques propri\'et\'es de l'homologie rigide qui d\'ecoulent directement de celles de la cohomologie rigide \`a support.

\begin{prop}
On a :
\begin{itemize}
\item  Soit $X$ une $k$-vari\'et\'e de dimension $n$ ; alors
$$H^{rig}_{i}(X) = 0 \hbox{ pour tout } i > 2n.$$
\item Si $X$ est une $k$-vari\'et\'e irr\'eductible de dimension $n$, l'espace $H^{rig}_{2n}(X)$ est canoniquement isomorphe \`a $K$. 
\item Si $X$ est une $k$-vari\'et\'e et $Z \hookrightarrow X$ un sous-sch\'ema ferm\'e, on a la suite exacte longue :
$$ \cdots \to H_{i+1}^{rig}(X-Z) \to H_{i}^{rig}(Z) \to H_{i}^{rig}(X) \to H_{i}^{rig}(X-Z) \to H_{i-1}^{rig}(Z) \to \cdots $$
\item Soit $K'$ une extension de $K$, d'anneau des entiers $\VCC'$ et de corps r\'esiduels $k'$. On note $X' = X\times_k k'$ le sch\'ema obtenu par changement de base. Il existe un isomorphisme canonique 
$$\varphi : H_{i}^{rig}(X'/K') \iso H_{i}^{rig}(X/K)\otimes_K K'.$$

\end{itemize}
\end{prop}

En utilisant la premi\`ere et la troisi\`eme assertion de la proposition on obtient :

\begin{coro}
Si $X$ est une $k$-vari\'et\'e et si on note $X_i$ ses composantes irr\'eductibles de dimension maximale $n$. On a :
$$H_{2n}^{rig}(X) = \bigoplus_{i} H_{2n}^{rig}(X_i).$$
\end{coro}

\subsection{Classe fondamentale}

Soit $Z$ une $k$-vari{\'e}t{\'e} int{\`e}gre de dimension $d$, on va d{\'e}finir sa classe fondamentale qui est un {\'e}l{\'e}ment de $H_{2d}^{rig}(Z)$. Pour cela on regarde le morphisme trace \cite{Berth7} : 
$$Tr_Z : H^{2d}_{c,rig}(Z) \to K.$$
Il d\'efinit une classe $\eta_Z \in H_{2d}^{rig}(Z)$.
Cette classe est fonctorielle par rapport aux immersions ouvertes car on a pour $U$ ouvert de $Z$ le diagramme commutatif suivant :
$$\xymatrix{ H^{2d}_{c,rig}(U) \ar[r]^-{Tr_U} \ar[d] & K \ar@{=}[d] \\  H^{2d}_{c,rig}(Z) \ar[r]^-{Tr_Z} & K. }$$

On a alors :

\begin{prop}
\label{1611}
Soit $K'$ une extension de $K$, d'anneau des entiers $\VCC'$ et de corps r\'esiduels $k'$. On note $X' = X\times_k k'$ le sch\'ema obtenu par changement de base et 
$$\varphi : H_{2n}^{rig}(X'/K') \to H_{2n}^{rig}(X/K)\otimes_K K'.$$
On a :
$$\varphi(\eta_{X'}) = \eta_{X}\otimes 1.$$
\end{prop}

\dem D'apr\`es la d\'efinition de la classe fondamentale, cette proposition se ram\`ene \`a montrer que le diagramme :
$$\xymatrix{ H^{2n}_{c,rig}(X'/K') \ar[rr]^{Tr_{X'}} & & K' \\ H^{2n}_{c,rig}(X/K)  \ar[u]^{\varphi} \ar[rr]_-{Tr_X} & & K  \ar[u] }$$
est commutatif, ce qui d\'ecoule de la d\'efinition du morphisme trace.
\findem

\begin{prop}
Soient $X$ et $Y$ deux $k$-vari\'et\'es de dimension $n$ et $m$ respectivement. Le morphisme de K\"{u}nneth :
$$H_{n+m}^{rig}(X\times_kY) \to H_n^{rig}(X) \otimes H_m^{rig}(Y)$$
envoie $\eta_{X\times Y}$ sur $\eta_X \otimes \eta_Y$.
\end{prop}

\dem Le r\'esultat d\'ecoule de la compatibilit\'e du morphisme de K\"{u}nneth au morphisme trace \cite{Berth7}.
\findem


\subsection{Cap-produit et formule de projection}
Nous allons d\'efinir le cap-produit et  d\'emontrer la formule de projection.

Pour toute immersion ferm{\'e}e $Y \hookrightarrow X$, nous d{\'e}finissons le {\it cap-produit} not{\'e} $\cap$ :
$$\begin{array}{ccc}
H^{rig}_i(X) \otimes H_{Y,rig}^j(X) & \to & H^{rig}_{i-j}(Y) \\
(\varphi, x) & \mapsto & (y \mapsto \varphi(x\cup y))
\end{array}$$
o{\`u} on a indentifi{\'e} $H^{rig}_i(X)$ (resp . $H^{rig}_{i-j}(Y)$) avec $H^i_{c,rig}(X)^{\vee}$ (resp. $H^{i-j}_{c,rig}(Y)^{\vee}$)  et not{\'e} $\cup$ le cup-produit \cite[2.2]{Berth7} :
$$ \cup : H^{j}_{Y, rig,}(X) \otimes H^{i-j}_{c,rig}(Y) \to H^i_{c,rig}(X).$$

Ce cup-produit est fonctoriel par rapport aux morphismes propres : pour tout diagramme cart\'esien 
$$\xymatrix{ Y' \ar@{^{(}->}[r]^{i'} \ar[d]_{f'} & X' \ar[d]^{f} \\ Y \ar@{^{(}->}[r]_{i} & X,}$$
o\`u $i$ et $i'$ sont des immersions ferm\'ees et $f$ et $f'$ des morphismes propres, et tous \'el\'ements $x \in H^j_{Y,rig}(X)$  et $y \in H^{i-j}_{c,rig}(Y)$ on a :
$$f^*(x \cup y) = f^*(x) \cup f^*(y).$$

On a alors 

\begin{prop}[Formule de projection]
Pour tout diagramme cart\'esien 
$$\xymatrix{ Y' \ar@{^{(}->}[r]^{i'} \ar[d]_{f'} & X' \ar[d]^{f} \\ Y \ar@{^{(}->}[r]_{i} & X,}$$
o\`u $i$ et $i'$ sont des immersions ferm\'ees et $f$ et $f'$ des morphismes propres, et tous \'el\'ements $x' \in H_i^{rig}(X')$ et $x \in H^j_{Y,rig}(X)$ on a :
$$f_*(x')\cap x = f'_*(x' \cap f^*(x)).$$
\end{prop}

\dem Notre {\'e}galit{\'e} se tenant dans $H^{rig}_{i-j}(Y) = H_{c,rig}^{i-j}(Y)^\vee,$ il suffit de montrer que pour tout $y \in H_{c,rig}^{i-j}(Y)$, on a
$$<f_*(x') \cap x, y> = <f'_*(x' \cap f^*(x)),y>$$
Or par d{\'e}finition du cap-produit, 
$$<f_*(x') \cap x, y> = <f_*(x'),x \cup y>.$$
De plus,

$$\begin{array}{ccc}
<f'_*(x' \cap f^*(x)),y> & = & <x' \cap f^*(x),f^*(y)> \\
  & = &  <x',f^*(x) \cup f^*(y)> \\
 & = & <x',f^*(x\cup y)> \\
 & = & <f'_*(x') , x \cup y >.
\end{array}$$

\findem

\subsection{Le cas des vari\'et\'es lisses}
Si on se donne $X$ une $k$-vari\'et\'e lisse et $Z$ un sous-sch\'ema ferm\'e int\`egre \'eventuellement singulier, il est plus habituel de regarder les groupes de cohomologie et non l'homologie. Les constructions pr\'ec\'edentes se r\'e\'ecrivent alors de la mani\`ere suivante.
On note $r = n-d$ la codimension de $Z$ dans $X$. La classe fondamentale de $Z$, encore not\'ee $\eta_Z$,  se construit donc comme un \'el\'ement de $H_{Z, rig}^{2r}(X)$. 

Dans ce cas, le cap-produit se r\'e\'ecrit en utilisant
$$H_i^{rig}(X) \iso H^{2n-i}_{rig}(X) \hbox{ et } H_{i-j}^{rig}(Z) \iso H^{2n-i+j}_{Z,rig}(X).$$
On obtient alors le cup-produit classique :
$$ H^{2n-i}_{rig}(X) \otimes H^j_{Z,rig}(X) \surfleche{\cup}  H^{2n-i+j}_{Z,rig}(X).$$
L'axiome de dualit\'e de Poincar\'e du th\'eor\`eme est alors ramen\'e \`a la dualit\'e de Poincar\'e \'enonc\'ee dans \cite{Berth7}.

Pour finir, la formule de projection s'exprime alors :
\begin{prop}[Formule de projection]
Pour tout diagramme cart\'esien 
$$\xymatrix{ Y' \ar@{^{(}->}[r]^{i'} \ar[d]_{f'} & X' \ar[d]^{f} \\ Y \ar@{^{(}->}[r]_{i} & X,}$$
o\`u $i$ et $i'$ sont des immersions ferm\'ees et $f$ et $f'$ des morphismes propres, et tous \'el\'ements $x' \in H^i_{rig}(X')$ et $x \in H^j_{Y,rig}(X)$ on a :
$$f_*(x')\cup x = f'_*(x' \cup f^*(x))$$
o\`u $f_* : H^i_{rig}(X') \to H^i_{rig}(X)$ est la transpos\'ee par les accouplements de Poincar\'e du morphisme $f^* : H^i_{c,rig}(X') \to  H^i_{c,rig}(X).$
\end{prop}

\subsection{Remarque sur la contravariance \'etale}
Les axiomes habituels des th\'eories de dualit\'e de Bloch-Ogus demandent une contravariance de l'homologie par rapport  aux morphismes \'etales et pas seulement aux  immersions ouvertes. De m\^emes, les fonctorialit\'e des suites exactes sont requises dans le cas des morphismes \'etales. Cependant nous n'utiliserons dans la suite que la contravariance par rapport aux immersions ouvertes.

\section{Classes de cycles}
Nous allons maintenant construire les classes de cycles \`a partir de la classe fondamentale.

\subsection{D\'efinition}
Commen\c{c}ons par rappeler les principales d\'efinitions.
Nous appelons groupe des cycles d'une $k$-vari{\'e}t{\'e} $X$ et notons
$Z(X)$ le groupe libre engendr{\'e} par les sous-sch{\'e}mas ferm{\'e}s
int{\`e}gres de $X$. On note $[T]$ la classe du sous-sch\'ema $T$ dans $Z(X)$. On gradue ce groupe par la dimension :
$$Z_\bullet(X) = \bigoplus_{k} Z_k(X),$$
o\`u $Z_k(X)$ est le sous-groupe engendr\'e par les sous-sch{\'e}mas ferm{\'e}s
int{\`e}gres de dimension $k$.

Rappelons qu'on associe alors {\`a} tout sous-sch{\'e}ma
ferm{\'e} $Z$ un cycle :
$$[Z]= \sum_{i=1}^n l(\OCC_{Z,t_i})[T_i^{red}]$$
o{\`u} les $T_i$ sont les composantes irr{\'e}ductibles de $Z$ et $t_i$
les points g{\'e}n{\'e}riques de ces derni{\`e}res.
Rappelons aussi les morphismes de fonctorialit{\'e} :

\begin{itemize}
\item
tout morphisme propre $f : X \to Y$ induit un morphisme
  $f_* : Z(X) \to Z(Y)$ d{\'e}fini sur les sous-sch{\'e}mas ferm{\'e}s
  int{\`e}gres de la mani{\`e}re suivante :
$$ f_{*}([Z]) = \left \{
\begin{array}{cc} 
0 & \hbox{si } \dim f(Z) < \dim Z \\
~[k(Z): k(f(Z))] [f(Z)] & \hbox{si } \dim f(Z) = \dim Z \\
\end{array} 
\right. $$
\item tout morphisme plat $f : X \to Y$ induit un morphisme $f^* : Z(Y)
  \to Z(X)$ d{\'e}fini sur les sous-sch{\'e}mas ferm{\'e}s
  int{\`e}gres de la mani{\`e}re suivante :
$$f^*([Z])=[Z \times_Y X].$$
\end{itemize}

Par lin\'earit\'e, on d\'efinit alors le morphisme {\it classe de cycle} :
$$\begin{array}{cccc}
\gamma :&  Z_\bullet(X) &\to& H_{\bullet}^{rig}(X) \\
   & \sum n_i [T_i] & \mapsto & \sum n_i (\alpha_i)_{*}\eta_{T_i},
\end{array}$$
o\`u $\alpha_i$ est le morphisme $T_i \hookrightarrow X$.

\subsection{Propri\'et\'es des classes de cycles}
Nous allons d\'emontrer que les classes de cycles sont fonctorielles et nous regarderons l'action du frobenius. Pour cela, nous allons devoir nous ramener au cas d'un diviseur lisse dans une vari\'et\'e lisse. L'ingr\'edient principal est le lemme suivant :

\begin{lemme}
\label{oter}
Soient $X$ une $k$-vari\'et\'e et $Z$ une sous-vari\'et\'e de codimension $r$. Il existe une extension finie $k'$ de $k$ telle que si on note $X'$ et $Z'$ les $k'$-vari\'et\'es d\'eduites de $X$ et $Z$ par extension des scalaires, il existe un ouvert $U'$ de $X'$ v\'erifiant, en notant $T'$ son compl\'ementaire :

\begin{itemize}
\item $Z' \cap T'$ est de codimension sup{\'e}rieure ou {\'e}gale {\`a} $r+1$ dans $X'$
\item $Z'_U = Z' \cap U'$ est lisse sur $k'$ 
\item Il existe $\ZCC$ et $\XCC$ deux $C'$-sch{\'e}mas
  affines et lisses et une
  immersion ferm{\'e}e $i : \ZCC \hookrightarrow \XCC$ se r{\'e}duisant
  sur $i_U : Z'_U \hookrightarrow U'$ o\`u $C'$ est un anneau de Cohen pour $k'$.
\end{itemize}
\end{lemme}

\dem

D'apr{\`e}s \cite[17.15.13]{EGAIV}.

\findem

\begin{prop}Soient $X$ et $X'$ deux $k$-vari\'et\'es.
Si $f: X \to X'$ est un morphisme propre et $x \in Z_{\bullet}(X)$, on a
$$\gamma(f_*(x)) = f_*(\gamma(x)).$$
\end{prop}

\dem 

Tout {\'e}tant lin{\'e}aire, on se ram{\`e}ne au cas d'un sous-sch{\'e}ma
ferm{\'e} int{\`e}gre $Z$. On note $Z'=f(Z)$, $\dim Z = d$, $\dim X = n$ et $\dim X' = n'$. On diff{\'e}rencie deux cas.
\begin{itemize}
\item Si $\dim Z' < d$, par d{\'e}finition $f_*(Z) = 0$. On est donc ramen{\'e} {\`a} montrer que $f_*(\gamma(Z)) = 0$. Le morphisme $f_{|Z} : Z \to Z'$ est propre, on a donc le diagramme commutatif suivant :
$$\xymatrix{ H_{2d}^{rig}(Z) \ar[r]^{i_*} \ar[d]^{(f_{|Z})_*} & H_{2d}^{rig}(X) \ar[d]_{f_*} \\
 H_{2d}^{rig}(Z') \ar[r]_{i'_*} & H_{2d}^{rig}(X').}$$
Or on sait que, par d{\'e}finition, $\gamma(Z)=i_*\eta_Z$. Donc on a 
$$f_*(\gamma(Z)) = i'_*((f_{|Z})_*(\eta_Z)).$$
Or $H_{2d}^{rig}(Z') = 0$. On a donc bien :
$$f_*(\gamma(Z)) = 0.$$
\item Si $\dim Z' = d$, il suffit gr\^ace \`a la proposition \ref{1611} de d\'emontrer cette \'egalit\'e apr\'es une extension finie des scalaires. On peut donc supposer, en utilisant le th\'eor\`eme de puret\'e et le lemme \ref{oter}, que $f$ est
un morphisme fini entre sch{\'e}mas affines et lisses. Comme $X$ est int\`egre, on sait qu'il existe un ouvert dense $U$ de $X$ au dessus duquel $f$ est plat. On note $U' = f^{-1}(U)$. Se restreindre \`a $U$ revient \`a \^oter  de $X$ et $X'$ des sous-sch\'emas ferm\'es de codimension sup\'erieure ou \'egale \`a $1$. On peut donc supposer que $f$ est fini et plat de degr\'e $\alpha=[k(X') : k(X)]$.  La compatibilit\'e des morphismes trace de \cite[2.3 et 1.4]{Berth7} se g\'en\'eralise aux morphismes finis et plats. On d\'eduit alors de \cite[3.6]{Berth3} le diagramme commutatif suivant : 
$$\xymatrix{ K \ar[r]^{\times \alpha} & K \\
H_{2d}^{rig}(X') \ar[r]_{f_*} \ar[u]_{Tr_{X'}} & H_{2d}^{rig}(X). \ar[u]_{Tr_{X}}},$$
ce qui implique la proposition.

\end{itemize}

\findem

On regarde $F : X \to X$ le morphisme de Frobenius absolu. Il induit 
$$\Phi : H^{rig}_i(X) \to H^{rig}_i(X).$$
On a alors 
\begin{prop}
Soit $Z$ une sous-vari\'et\'e de codimension $r$ dans $X$. On a 
$$\Phi(\gamma(Z)) = p^r \gamma(Z).$$
\end{prop}
\dem On peut par lin\'earit\'e se ramener au cas d'un sous sch\'ema ferm\'e int\`egre $Z$. Si on note $n$ la dimension de $X$. On va montrer la propri\'et\'e pour la classe fondamentale
$\eta_Z \in H^{rig}_{2(n-r)}(Z)$.

De plus, gr\^ace \`a la proposition \ref{1611}, on peut utiliser le lemme \ref{oter} pour se ramener au cas d'un sous-sch\'ema lisse dans une vari\'et\'e lisse. Dans ce cas, $\eta_Z$ se construit naturellement comme \'el\'ement de $H^{2r}_{Z,rig}(X)$. On peut, de plus, supposer que l'immersion ferm\'ee $Z \hookrightarrow X$ est la fibre sp\'eciale d'une immersion ferm\'e de $C$-sch\'ema $\ZCC \hookrightarrow \XCC$ d\'efinie par $r$ sections globales sur $\XCC$. 

\begin{lemme}
Avec les notations pr\'ec\'edentes, on note 
$$G^{rig}_{Z/X} : H^0_{rig}(Z) \to H_{Z,rig}^{2r}(X)$$
le morphisme de Gysin \cite{Berth7}.
On a alors
$$\eta_Z = G^{rig}_{Z/X}(1).$$
\label{finito}
\end{lemme}

\dem cela d\'ecoule directement de la compatibilit\'e entre le morphisme de Gysin et l'accouplement de dualit\'e.
\findem

Dans ce cas, le th\'eor\`eme 2.4 de \cite{Chia}, nous dit que l'on a le diagramme commutatif suivant :
$$\xymatrix{H^0_{Z,rig}(X) \ar[r]^{G_{Z/X}} \ar[d]_{p^r\Phi} & H^{2r}_{Z,rig}(X) \ar[d]^{\Phi} \\ H^0_{Z,rig}(X) \ar[r]^{G_{Z/X}}  & H^{2r}_{Z,rig}(X).}$$ 
Notre proposition en d\'ecoule.

\findem

\section{Compatibilit\'e \`a l'\'equivalence rationnelle}
Nous allons montrer que les classes de cycles d\'efinies ci-dessus, sont compatibles \`a l'\'equivalence rationnelle. Pour ce faire, nous aurons besoin de la trivialit\'e de la classe d'un diviseur principal ({\it principal triviality} dans \cite{Blo-Ogu, Jann}).

\subsection{Classe d'un diviseur}

Nous allons rappeler les d\'efinitions et r\'esultats principaux sur les diviseurs. On trouvera un expos\'e plus pr\'ecis sur le sujet dans \cite[II.6]{Hart1} et \cite[21]{EGAIV}.

Soit $X$ une $k$-vari\'et\'e, on note $\MCC$ le faisceau associ\'e au pr\'efaisceau $\FCC$ qui est d\'efini sur tout ouvert affine $U=\spec{A}$ de $X$ par :
$$\FCC(U) := S^{-1}A$$
o\`u $S$ est le syst\`eme multiplicatif des \'el\'ements $f \in A$ tel que pour tout $x \in U$, $f$ est non diviseur de z\'ero dans l'anneau local $\OCC_{U,x}$.
Le faisceau $\MCC$, qui peut \^etre d\'efini sur tout sch\'ema, est la g\'en\'eralisation du corps des fonctions sur un sch\'ema int\`egre.

\begin{Defi}
Avec les notations pr\'ec\'edentes, on appelle diviseur de Cartier une section globale du faisceau quotient $\MCC^*/\OCC^*$. De plus, un diviseur de Cartier sera dit principal s'il appartient \`a l'image de l'application
$$\Gamma(X, \MCC^*) \to \Gamma(X, \MCC^*/\OCC^*).$$
\end{Defi}

Un diviseur de Cartier peut donc \^etre repr\'esent\'e par un recouvrement $\{U_i\}$ de $X$ et la donn\'ee pour tout $i$ de fonctions m\'eromorphes $f_i \in \Gamma(U_i, \MCC^*)$, telles que pour tout couple $i, j$, on ait 
$$\frac{f_i}{f_j} \in \Gamma(U_i \cap U_j, \OCC_X^*).$$

On associe \`a tout diviseur de Cartier $D$ un faisceau inversible $\LCC(D)$ en prenant le sous-faisceau de $\MCC$ engendr\'e par $f_i^{-1}$ sur $U_i$.

\begin{Defi}
Un diviseur de Cartier $D$ sur $X$ sera dit effectif s'il peut \^etre repr\'esent\'e par la donn\'ee $\{(U_i,f_i)\}$ o\`u les $U_i$ sont des ouverts qui recouvrent $X$ et $f_i \in \Gamma(U_i, \OCC_X)$. Dans ce cas on lui associe un sous-sch\'ema de codimension $1$, $Y(D)$ qui est le sous-sch\'ema d\'efini par le faisceau d'id\'eaux $\ICC(D)$ localement engendr\'e par $f_i$.
\end{Defi}

Soit $D$ un diviseur de Cartier sur $X$, on notera $|D|$ le support de $D$ c'est \`a dire le sous-sch\'ema ferm\'e des points $x \in X$ tel que $D_x \neq 0$.

Si $X$ est une vari\'et\'e lisse, le groupe des diviseurs de Cartier est isomorphe au groupe des diviseurs de Weil et les diviseurs de Cartier effectifs correspondent aux diviseurs de Weil effectifs. Nous noterons $[D]$ le diviseur de Weil associ\'e au diviseur de Cartier $D$.

On consid{\`e}re maintenant $X$ une $k$-vari{\'e}t{\'e} lisse et $D$ un diviseur ; on regardera $\LCC(D)$ le faisceau inversible qui lui est canoniquement associ{\'e}. On sait lui associer sa premi\`ere classe de Chern \cite{Pet1}.

Si on note $U$ le compl\'ementaire de $D$ dans $X$, on voit que l'image par l'application canonique
$$H^2_{rig}(X) \to H^2_{rig}(U)$$
de la premi\`ere classe de Chern $c_1(\LCC(D))$ est nulle car $\LCC(D)$ est trivialis\'e sur $U$.

Notre but est de relever la premi{\`e}re classe de Chern de $\LCC(D)$ en un  \'el\'ement $c_{1,D}$ dans le groupe de cohomologie rigide {\`a} support $H^2_{|D|, rig}(X).$ Nous appelerons cet \'el\'ement {\it classe du diviseur} $D$.

Notre diviseur de Cartier, s'\'ecrit de mani\`ere unique 
$$D = D^+ - D^-$$
avec $D^+$ et $D^-$ des diviseurs effectifs dont les supports ne contiennent pas de composantes communes. Cette derni\`ere condition nous assure que 
$$H^2_{|D|,rig}(X) \iso H^2_{|D^+|,rig}(X) \oplus  H^2_{|D^-|,rig}(X).$$
On posera alors 
$$c_{1,D} = c_{1,D^+} - c_{1,D^-}.$$

On est donc ramen\'e au cas des diviseurs effectifs.

Rappelons pour commencer la construction de la cohomologie rigide \cite{Berth3,Berth7}. Soient $X$ une vari\'et\'e et $Z$ un sous-sch\'ema ferm\'e. On choisit une compactification $\ov{X}$ de $X$ et on suppose qu'il existe une immersion ferm\'ee de $\ov{X}$ dans un $C$-sch\'ema formel $\YCC$ lisse au voisinage de $X$. On note $U = X - Z$, $j_ X \hookrightarrow \ov{X}$, $j_U : U \hookrightarrow \ov{X}$ et $\ICC$ l'id\'eal de $\OCC_{\YCC}$ d\'efinissant $\ov{X}$. En notant $j_X^\dag$ et $j_U^\dag$, les foncteurs d\'efinis dans \cite[1.2]{Berth3}, les groupes de cohomologie \`a support se calculent de la mani\`ere suivante :
$$H^i_{Z,rig}(X) := \HM^i(]X[, (j_X^\dag\Omega_{]X[}^\star \to j_U^\dag\Omega_{]X[}^\star)_s).$$
Dans la formule ci-dessus, $(j_X^\dag\Omega_{]X[}^\star \to j_U^\dag\Omega_{]X[}^\star)_s$ d\'esigne le complexe simple associ\'e au complexe double o\`u le terme de bidegr\'ee $(0,0)$ est $j_X^\dag\OCC_{]X[}$. La diff\'erentielle $d : j_X^\dag\Omega^i_{]X[} \oplus j_U^\dag\Omega^{i-1}_{]X[} \to j_X^\dag \Omega^{i+1}_{]X[} \oplus j_U^\dag\Omega^{i}_{]X[}$ est donn\'ee par
$$\left(\begin{array}{cc} d_X & 0 \\ u & -d_U\end{array} \right),$$
o\`u $u$ est la restriction et $d : j_X^\dag \Omega^{i-1}_{]X[} \to j_X^\dag \Omega^{i}_{]X[}$ (resp. $d_U : j_U^\dag \Omega^{i-1}_{]X[} \to j_U^\dag \Omega^{i}_{]X[}$ ) est obtenue en appliquant le foncteur $j_X^\dag$ (resp. $j_U^\dag$ ) \`a $d_X$.

On va maintenant \'etudier en premier lieu le cas des vari\'et\'es propres. On ram\`enera alors le cas g\'en\'eral \`a un \'enonc\'e d'ind\'ependance de prolongements.

Soient $X$ une vari\'et\'e propre int\`egre mais {\bf non n\'ecessairement lisse} et $D$ un diviseur de Cartier effectif. On note $j_U : U = X - |D| \hookrightarrow X.$ 

On va mener une construction similaire \`a \cite[2.1]{Pet1}.

On se donne un repr\'esentant de $D$ : $\{(X_i), h_i\}$ et un plongement de $X$ dans un $C$-sch\'ema formel $\YCC$ lisse au voisinage de $X$ d\'efini par un id\'eal $\ICC$. On sait qu'il existe un recouvrement fini affine $\UG=(\YCC_i)_{i\in I}$ de $\YCC$ dont la restriction \`a $X$ est un raffinement du recouvrement $(X_i)$. On peut donc supposer que $X_i = \YCC_i \cap X$.

On note $\YCC_i = \spf{\ACC_i}$ les ouverts du recouvrements $\UG$ et $X_i = \spec{A_i}$ les ouverts induits sur $X$. On notera de plus $U_i$ l'ouvert $D(h_i)$ de $X_i$. Le faisceau $\LCC(D)$ est trivialis\'e sur ce recouvrement. Par d\'efinition, le cocycle
$$u_{ij} = \frac{h_i}{h_j}$$
est un repr\'esentant de $\LCC(D)$.
On se donne alors un rel\`evement $\tilde{h}_i \in \ACC_i$ de $h_i$ et un rel\`evement $\tilde{u}_{ij}$ de $u_{ij}$ dans $\ACC_{ij}$. On a dans $\ACC_{ij}$ :
\begin{equation}
\label{argen}
\tilde{h}_j.\tilde{u}_{ij} = \tilde{h}_{i} + \alpha_{ij}
\end{equation}

avec $\alpha_{ij} \in \Gamma(U_{ij}, \ICC).$ On regarde $\tilde{h}_i$ et $\widetilde{u}_{ij}$ comme des fonctions analytiques sur $]X_i[$ et $]X_{ij}[$ respectivement. 

Il existe alors (lemme 2.1.2 de \cite{Pet1}) un voisinage strict $V_{ij}$ de $]U_{ij}[$ dans $]X_{ij}[$ sur lequel la fonction $\tilde{h}_i$ est inversible. De plus on peut supposer que pour tout $x \in V_{ij}$ on a, en posant $\mu_{ij} = \alpha_{ij}/\tilde{h}_i$ :
$$|\mu_{ij}(x)| < 1.$$

On peut alors regarder l'\'equation \ref{argen} comme une \'egalit\'e de fonctions analytiques sur $V_{ij}$. On peut \'ecrire notre \'equation sous la forme :
$$\tilde{h}_j.\tilde{u}_{ij} = \tilde{h}_{i}.(1+\mu_{ij}).$$

Par suite, on sait \cite{Pet1} que $\tilde{u}_{ij}$ est inversible. On d\'efinit
$$c_{1,D} \in C^2(\UG_K,(\Omega_{]X[}^\star \to j_U^\dag\Omega_{]X[}^\star)_s) = C^2(\UG_K,\OCC_{]X[}) \oplus C^1(\UG_K,\Omega^1_{]X[} \oplus j_U^\dag\OCC_{]X[}) \oplus C^0(\UG_K,\Omega^2_{]X[} \oplus j_U^\dag\Omega^1_{]X[})$$
en posant :
$$(c_{1,D})_{ijk} := - \log (\tilde{u}_{ij}\tilde{u}_{ik}^{-1}\tilde{u}_{jk}) \in  C^2(\UG_K,\OCC_{]X[}),$$
$$(c_{1,D})_{ij} := \frac{d\tilde{u}_{ij}}{\tilde{u}_{ij}}+ \log(1+\mu_{ij}) \in C^1(\UG_K,\Omega^1_{]X[} \oplus j_U^\dag\OCC_{]X[})$$
et 
$$(c_{1,D})_i := -\dlog{\tilde{h}_i} \in  C^0(\UG_K,\Omega^2_{]X[} \oplus j_U^\dag\Omega^1_{]X[}).$$

Avec cette d\'efinition, il est imm\'ediat en reprenant  les calculs de \cite[2]{Pet1} de v\'erifier que l'on construit bien une classe ne d\'ependant pas de nos choix.

\begin{prop}
Soient $f : X' \to X$ un morphisme de vari\'et\'es propres int\`egres et $D$ un diviseur de Cartier effectif sur $X$ tel que $f^*(D)$ existe, alors :
$$f^*c_{1,D} = c_{1,f^*D}.$$
\end{prop}
\dem Il suffit de reprendre la construction pr\'ec\'edente pour v\'erifier qu'elle est fonctorielle.
\findem

On va maintenant revenir au cas g\'en\'eral. Soit $X$ une vari\'et\'e lisse mais non n\'ec\'essairement propre. On peut, en travaillant composante par composante supposer que $X$ est int\`egre. 
\begin{lemme}
Avec les notations pr\'ec\'edentes, il existe une compactification $X \hookrightarrow \ov{X}$ et un diviseur de Cartier effectif $\ov{D}$ sur $\ov{X}$ dont la restriction \`a $X$ est $D$.
\end{lemme}
\dem
Il existe une compactification $X \hookrightarrow \ov{X}$, o\`u $\ov{X}$ est int\`egre et $\ov{X} - X$ est un diviseur \cite[4.1.1]{Pet1}. On consid\`ere alors $\ov{|D|}$ l'adh\'erence sch\'ematique de $|D|$ dans $\ov{X}$. On \'eclate $\ov{|D|}$ dans $\ov{X}$. On voit alors que la vari\'et\'e obtenue, $\ov{X}'$ est encore propre et que l'image inverse de $\ov{X} - X$ est un diviseur. L'image inverse de $\ov{|D|}$ est alors un diviseur de Cartier not\'e $\ov{D}$. De plus, si on restreint l'\'eclatement pr\'ec\'edent \`a $X$, on ne modifie rien puisqu'on \'eclate un diviseur de Cartier. On peut donc choisir $\ov{X}'$ et $\ov{D}$.
\findem

A partir de l\`a on veut d\'efinir, comme pour les classes de Chern, la classe d'un diviseur par
$$c_{1,D} = j^*c_{1,\ov{D}},$$
o\`u $j : X \hookrightarrow \ov{X}$.

\begin{prop}
Soient $\ov{X}$ une vari\'et\'e propre et int\`egre, $X$ un ouvert  et $\ov{D}_1$, $\ov{D}_2$ deux diviseurs de Cartier effectifs tels que 
$$(\ov{D}_1)_{|X} = (\ov{D}_2)_{|X} = D,$$
alors on a dans $H^2_{|D|,rig}(X)$ :
$$j^*c_{1,\ov{D}_1} = j^*c_{1,\ov{D}_2},$$
o\`u $j : X \hookrightarrow \ov{X}.$
\end{prop}

\dem La d\'emonstration est similaire \`a celle du th\'eor\`eme 4.2.1 de \cite{Pet1}. Nous allons nous contenter de mettre en avant les diff\'erences. Nous reprendrons les m\^emes notations. On se donne un plongement de $\ov{X}$ dans un $C$-sch\'ema formel $\YCC$ lisse aux points de $\ov{X}$. On choisit un recouvrement affine fini $\UG = (\YCC_i)_{i\in I}$ dont la restriction \`a $\ov{X}$ trivialise les deux diviseurs de Cartier. On note $h_i^1$ (resp. $h_i^2$) l'\'equation locale d\'efinissant $\ov{D}_1$ (resp. $\ov{D}_2$) dans $\ov{X}_i$. Par hypoth\`eses, on sait qu'il existe pour tout $i \in I$, un \'el\'ement $\theta_i \in (A_i)_{f_i}^*$ tel que l'on ait dans $(A_i)_{f_i}^*$ :
$$h_i^2=h_i^1.\theta_i.$$
Cette notation est coh\'erente avec celle de \cite{Pet1}.

En proc\'edant comme pr\'ec\'edemment : on multiplie par des puissances de $f_i$ puis on rel\`eve, on obtient l'existence de $k \in \NM$ et $\kappa_i \in I_i$ tels que
\begin{equation}
\label{gal}
f_i^k\tilde{h}^2_i = \tilde{h}^1_i.\tilde{\alpha}_i + \kappa_i,
\end{equation}
o\`u $\tilde{\alpha}_i$ est un rel\`evement de $\alpha:= \theta_i.f_i^k$.
On note $U$ le compl\'ementaire de $D$ dans $X$ et $U_i$ l'intersection de $U$ avec $\YCC_i$.
On se restreint alors \`a des voisinages stricts $V_i$ de $]U_i[$ dans $]\ov{X}_i[$. On peut donc supposer que $\tilde{h}^1_i$ et $\tilde{h}^2_i$ sont inversibles. 
La restriction de l'\'equation \ref{gal} \`a ce voisinage strict peut donc s'\'ecrire :
$$ \tilde{h}^2_i = \tilde{h}^1_i.\tilde{\theta}_i\left( 1 + \frac{\kappa_i}{f_i^k}(\tilde{h}^1_i.\tilde{\theta}_i)^{-1}\right),$$
o\`u on a pos\'e $\tilde{\theta}_i = \tilde{\alpha}_i/f_i^k.$
On montre que, quitte \`a se restreindre \`a des voisinages stricts $V_i$ de $]U_i[$, pour tout $x \in V_i$ :
$$\left|\frac{\kappa_i}{f_i^k}(\tilde{h}^1_i.\tilde{\theta}_i)^{-1}(x)\right| < 1.$$

On pose alors 
$$\omega_i := -\log(1+\frac{\kappa_i}{f_i^k}(\tilde{h}^1_i.\tilde{\theta}_i)^{-1}) \in C^0(\UG_K, j_U^\dag \OCC_{]\ov{X}[}).$$
On note
$$\zeta_{i} := \dlog{\tilde{\theta_i}} + \omega_i \in  C^0(\UG_K, j_X^\dag \Omega^1_{]\ov{X}[}\oplus j_U^\dag \OCC_{]\ov{X}[}),$$
et 
$$\zeta_{ij} \in  C^1(\UG_K, j^\dag\OCC_{]\ov{X}[}),$$
l'\'el\'ement d\'efini dans la d\'emonstration du th\'eor\`eme 4.2.1 de \cite{Pet1}. $\zeta_i$ et $\zeta_{ij}$ d\'efinissent un \'el\'ement 
$$\zeta \in C^1(\UG_K, (j_X^\dag \Omega^\star_{]\ov{X}[} \to j_U^\dag \Omega^\star_{]\ov{X}[})_s).$$

On a 
$$j^*c_{1,\ov{D}_1} - j^*c_{1,\ov{D}_2} = \Delta(\zeta).$$

\findem

On a donc associ\'e \`a tout diviseur de Cartier $D$ sur $X$ une classe 
$$c_{1,D} \in H^2_{|D|,rig}(X)$$

\begin{prop}
\label{galois}
Soient $f : X' \to X$ un morphisme de vari\'et\'es lisses  et $D$ un diviseur de Cartier sur $X$ tel que $f^*(D)$ existe, alors :
$$f^*c_{1,D} = c_{1,f^*D}.$$
\end{prop}





\subsection{Trivialit\'e de la classe d'un diviseur principal}
On se donne une vari\'et\'e lisse $X$ de dimension $n$ et $D$ un diviseur sur $X$. Soit $[D] = \sum n_i Z_i$ le diviseur de Weil associ\'e. On peut d\'efinir la classe fondamentale de $[D]$ par lin\'earit\'e :
$$\eta_{[D]} := \sum n_i \eta_{Z_i} \in \bigoplus H^{2}_{Z_i,rig}(X) = H^2_{|D|,rig}(X).$$

Nous allons comparer la classe de
fondamentale $\eta_{[D]}$ avec la classe $c_{1,D}$ du diviseur $D$.

Nous allons avoir besoin d'espaces analytiques. Nous allons donc
reprendre la d{\'e}finition \cite[0.3]{Berth5}

Nous allons d\'efinir la structure analytique pour les vari\'et\'es affines, le cas g\'en\'eral s'obtenant par recollement. 

Soit $Y$ une $K$-vari{\'e}t{\'e} affine. On pose $Y=\spec{A}$ et on se
donne une pr{\'e}sentation
$$A = K[T_1, \ldots, T_q]/(f_1, \ldots f_r).$$

On pose alors $T^{(m)}_i= \pi^mT_i$, $f^{(m)}_j=f_j(\pi^{-m}T_1^{(m)},
\ldots,\pi^{-m}T_q^{(m)})$ et 
$$\widehat{A}_m := K\{T^{(m)}_1, \ldots, T^{(m)}_q\}/(f^{(m)}_1, \ldots,  f^{(m)}_r).$$

On note $Y^f$ l'ensemble des points ferm{\'e}s de $Y$ et $Y_m =
\{x\in Y^f | |T_i(x)| \leq |\pi|^{-m}\}$. On a alors une bijection
$\spm{\widehat{A}_m} \to Y_m$ qui munit ce dernier d'une structure
d'espace analytique. On utilise pour finir que $Y^f = \cup_{m} X_m$
pour munir $Y^f$ d'une structure d'espace analytique. Cette structure
ne d{\'e}pend ni du choix de $\pi$ ni de celui de la pr\'esentation. On notera $Y^{an}$ la vari{\'e}t{\'e} analytique ainsi d{\'e}finie.

\rqe il faut faire attention que si la structure de vari\'et\'e analytique ainsi d\'efinie ne d\'epend pas du choix de la pr\'esentation les ouverts $Y_m$ eux en d\'ependent.

Il existe aussi un foncteur {\it faisceau analytique} :

$$\begin{array}{ccc} {\rm Mod}(Y) & \to & {\rm Mod}(Y^{an}) \\
    \FCC & \mapsto & \FCC^{an} 
\end{array}$$
o{\`u} ${\rm Mod}(Y)$ (resp. ${\rm Mod}(Y^{an})$) d{\'e}signe la cat{\'e}gorie des
faisceaux en $\OCC_Y$ (resp. $\OCC_{Y^{an}}$) modules sur $Y$ (resp. $Y^{an}$).
On peut d{\'e}crire ce foncteur de la mani{\`e}re suivante. Pour tout
ouvert $U$ de $Y$, $U^{an}$ est un ouvert de $Y^{an}$ et on a :
$$\FCC^{an}(U^{an}) = \FCC(U) \otimes_{\OCC_Y(U)} \OCC_{Y^{an}}(U^{an}).$$

\begin{prop}
Le foncteur faisceau analytique associ{\'e} est exact.
\end{prop}

\dem Cela d{\'e}coule du fait que $\OCC_{Y^{an}}$ est plat sur $\OCC_Y$.

\findem

On va d\'emontrer alors :

\begin{prop}
\label{prop331}
Soient $X$ une vari\'et\'e lisse et $D$ un diviseur sur $X$. On a dans $H^2_{|D|,rig}(X)$, l'\'egalit\'e :

$$\eta_{[D]} = c_{1,D}.$$

\end{prop}
\dem
Par lin\'earit\'e, on peut supposer que $D$ est le diviseur associ\'e \`a un sous-sch\'ema ferm\'e int\`egre $Z$. De plus la classe fondamentale et la classe d'un diviseur \'etant compatible \`a l'extension des scalaires, il nous suffit de d\'emontrer cette \'egalit\'e apr\'es une extension finie du corps de base.

On peut donc utilise le lemme \ref{oter}. Soit $U$ un ouvert v\'erifiant les hypoth\`eses du lemme, on a :

$$c_{1,D} = c_{1,D_U} \hbox{ et } \eta_{[D]} = \eta_{[D_U]}$$
via l'isomorphisme :
$$H^2_{|D|, rig}(X) \iso H^2_{|D_U|, rig}(U).$$

On peut donc supposer que $Z$ et $X$ sont affines et lisses et qu'il existe deux $C$-sch{\'e}mas $\ZCC$ et $\XCC$
affines  et lisses ainsi qu'une immersion ferm{\'e}e de $\ZCC$ dans $\XCC$
telle que l'on obtienne l'inclusion de $Z$ dans $X$ en passant aux
fibres sp{\'e}ciales. Quitte {\`a} localiser encore, on peut de plus supposer qu'il existe $t \in \Gamma(\XCC, \OCC_\XCC)$ tel que $\ZCC = V(t).$
On note respectivement $Z_K$ et $X_K$
les fibres g{\'e}n{\'e}riques de $\ZCC$ et $\XCC$.

On sait alors gr\^ace au lemme \ref{finito} que la classe fondamentale est l'image de $1$ par le morphisme de Gysin. On va exprimer explicitement le morphisme de Gysin alg{\'e}brique de
l'immersion $\ZCC \hookrightarrow \XCC$. On sait d'apr{\`e}s \cite[VI.3]{Berth1} qu'il existe un
morphisme de Gysin en cohomologie de De Rham :
$$G_{\ZCC/\XCC} : \Omega^\bullet_{\ZCC} \to \Omega^\bullet_{\XCC}[2].$$
On note alors $\HCC_{\ZCC}^1(\Omega^\bullet_{\XCC})$ le complexe :
$$0 \to \HCC_{\ZCC}^1(\OCC_{\XCC}) \to \HCC_{\ZCC}^1(\Omega^1_{\XCC}) \to \cdots \to
\HCC_{\ZCC}^1(\Omega^i_{\XCC}) \to \ldots $$ o{\`u} $\HCC_{\ZCC}^1$ d{\'e}signe le premier
groupe de cohomologie {\`a} support dans $\ZCC$. Le morphisme de Gysin se factorise
alors :
$$G_{\ZCC/\XCC} : \Omega^\bullet_{\ZCC} \to \HCC_{\ZCC}^1(\Omega^\bullet_{\XCC})[1] \to \Omega^\bullet_{\XCC}[2].$$

Nous allons expliciter ce morphisme en tant que morphisme de
complexes. 

Pour tout faisceau de $\OCC_{\XCC}$-modules $\FCC$ plat
sur $\OCC_{\XCC}$, on a la r{\'e}solution $\MCC(\FCC)$ suivante
\cite[2.2.2]{Hart3} :
\begin{equation}
\label{resol}
\FCC \to \MCC^0(\FCC)=i_*(i^*\FCC) \to \MCC^1(\FCC)=\HCC^1_{Z}\FCC
\end{equation}
o{\`u} $i : \XCC-\ZCC \hookrightarrow \XCC$ est l'immersion ouverte canonique.

Le morphisme de Gysin peut alors {\^e}tre vu comme le morphisme de
complexes suivant :

$$\Omega^\bullet_{\ZCC} \surfleche{a} \HCC^1_\ZCC(\OCC_\XCC)\otimes \Omega^\bullet_\XCC[1] \surfleche{b} \MCC(\Omega^\bullet_\XCC)[2] \liso
\Omega^\bullet_\XCC[2].$$

Le morphisme $a$  s'explicite de la mani{\`e}re suivante. Pour tout ouvert $W$,
 tout $\alpha$ et tout $\omega \in \Gamma(W,
 \Omega^\alpha_{\ZCC})$ on pose :
$$a(\omega) = \frac{\widetilde{\omega}}{t}\wedge dt \in
\Gamma(W, \HCC^1_{\ZCC}(\OCC_\XCC)\otimes\Omega^{\alpha+1}_{\XCC}) \subset
\Gamma(W, \MCC(\Omega_\XCC^\bullet)^{\alpha+2}).$$
o{\`u} $\widetilde{\omega}$ est un rel{\`e}vement de $\omega$ dans $\Gamma(W, \Omega^\alpha_{\XCC}).$ 

On consid{\`e}re le morphisme induit sur la fibre g{\'e}n{\'e}rique. On passe ensuite aux vari{\'e}t{\'e}s analytiques associ{\'e}es. On obtient alors le morphisme dans la cat{\'e}gorie des complexes de faisceaux en $K$-vectoriels sur $X_K^{an}$ :
$$G^{an}_{Z_K/X_K} : \Omega_{Z_K^{an}}^\bullet \to \HCC^1_{Z_K}(\OCC_{X_K})^{an} \otimes \Omega^\bullet_{X_K^{an}}[1] \to \MCC(\Omega^\bullet_{X_K^{an}})[2]$$
o{\`u} la r{\'e}solution $\MCC$ est obtenue en appliquant le foncteur exact {\it faisceau analytique associ{\'e}} {\`a} \ref{resol}.

On se donne maintenant $(\ov{\XCC},j)$ une compactification de $\XCC$ sur $C
$ et on note $\ov{\ZCC}$ l'adh{\'e}rence sch{\'e}matique de $\ZCC$ dans $\ov{\XCC}$. On note alors $\ov{X}$ (resp. $\ov{Z}$), $\ov{X}_K$ (resp. $\ov{Z}_K$)  et $\widehat{\ov{\XCC}}$ (resp. $\widehat{\ov{\ZCC}}$) la fibre sp{\'e}ciale, la fibre g{\'e}n{\'e}rique et le compl{\'e}t{\'e} formel de $\ov{\XCC}$ (resp. $\ov{\ZCC}$). Comme $\ov{\XCC}$ et $\ov{\ZCC}$ sont propres sur $C$ on a :
$$\widehat{\ov{\XCC}}_K = \ov{X}_K^{an} \hbox{ et } \widehat{\ov{\ZCC}}_K = \ov{Z}_K^{an}.$$

Le morphisme de Gysin rigide est obtenu en appliquant le foncteur $\Gamma^\dag_{]\ov{Z}[}j^\dag$ \cite[5.4]{Berth3} \`a $G^{an}_{Z_K/X_K}$.

On obtient un morphisme de complexes de faisceaux de $K$-vectoriels sur $\ov{X}_K^{an}$
$$\Gamma^\dag_{]\ov{Z}[}j^\dag\Omega^\bullet_{\ov{Z}_K^{an}} \to \Gamma^\dag_{]\ov{Z}[}j^\dag(\HCC^1_{\ov{Z}_K}(\OCC_{\ov{X}_K})^{an} \otimes \Omega^\bullet_{\ov{X}_K^{an}}[1]) \to \Gamma^\dag_{]\ov{Z}[}j^\dag(\MCC(\Omega^\bullet_{\ov{X}_K^{an}})[2]).$$
Or, les foncteurs $j^\dag$ et $\Gamma^\dag_{]\ov{Z}[}$ sont exacts, donc $\Gamma^\dag_{]\ov{Z}[}j^\dag(\MCC(\Omega^\bullet_{\ov{X}_K^{an}}))$ est quasi-isomorphe {\`a} $\Gamma^\dag_{]\ov{Z}[}j^\dag(\Omega^\bullet_{\ov{X}_K^{an}}).$

On consid{\`e}re pour finir $\zeta \in C^0(\UG, \Gamma^\dag_{]\ov{Z}[}j^\dag
\MCC^0\Omega^1_{X_K^{an}}) \subset C^1(\UG, \Gamma^\dag_{]\ov{Z}[}j^\dag
\MCC^\bullet\Omega^\bullet_{X_K^{an}})$ d{\'e}fini par :
$$\zeta_i = \frac{dt_i}{t_i}.$$

Un calcul montre que, si on note $d$ et $d'$ les d{\'e}rivations du
bicomplexe h{\'e}rit{\'e}es respectivement de la d{\'e}rivation de
$C^\bullet$ et de $\MCC^\bullet$, on a :
$$d\zeta = c_{1,Z} \hbox{ et } d'\zeta = G_{Z/X}^{rig}(1) = \eta_Z.$$
\findem

\begin{coro}
Soit $i : |D| \hookrightarrow X$ un diviseur principal lisse dans une vari\'et\'e lisse. On a $$i_*\eta_{[D]} = 0.$$
\end{coro}
\dem On a 
$$i_*\eta_{[D]} = i_*c_{1,D} = c_1(\LCC(D)).$$
Or $D$ est principal donc $\LCC(D)$ est trivialisable sur $X$. 
\findem

\subsection{Passage \`a l'\'equivalence rationnelle}

Nous allons commencer par rappeler la d\'efinition et les propri\'et\'es de l'homologie de Chow qui se trouvent dans \cite{Fult1, Fult2}.

Un cycle $z \in Z_{\bullet}(X)$ est dit rationnellement \'equivalent \`a z\'ero \cite{Fult2}, s'il existe un morphisme propre $\pi : Y \to X$ et un diviseur de Cartier principal $D$ sur $Y$ tel que 
$$z = \pi_*([D]).$$

Les cycles rationnellement \'equivalents \`a z\'ero forment un sous-groupe gradu\'e de $Z_{\bullet}(X)$ et le groupe quotient est appel\'e {\it homologie de Chow} et not\'e $A_\bullet(X)$. 

\begin{prop}
Les morphismes de fonctorialit\'e ci-dessus passent au quotient par l'\'equivalence rationnelle. Pr\'ecis\'ement, tout morphisme propre $f : X \to Y$ d\'efinit un morphisme $f_* : A_\bullet(X) \to A_\bullet(Y)$ et tout morphisme plat $f : X \to Y$  un morphisme $f^* : A_\bullet(Y) \to A_\bullet(X)$.
\end{prop}

\begin{prop}
Soient $X$ une $k$-vari\'et\'e et $z \in Z(X)$ un cycle rationnellement \'equivalent \`a z\'ero. On a 
$$\gamma(z) = 0.$$
\end{prop}

Avant de donner la d\'emonstration de la propri\'et\'e, \'enon\c{c}ons un corollaire.

\begin{coro}
Soit $X$ une $k$-vari\'et\'e, l'application classe de cycle passe au quotient par l'\'equivalence rationnelle et d\'efinit de la sorte :
$$\gamma : A_\bullet(X) \to H_{\bullet}^{rig}(X).$$
\end{coro}

\dem Il suffit de montrer que, si on se donne sur $X$ un diviseur de Cartier principal $D$, on a $\gamma([D]) = 0.$ L\`a encore grace \`a la compatibilit\'e des classes de cycles aux extensions finies du corps de base, nous allons pouvoir utiliser des techniques d'excisions.

Commen\c{c}ons par montrer que l'on peut supposer que $X$ est normale. En effet soit $\pi : X' \to X$ la normalisation de $X$. On sait \cite[1.5]{Fult2} que :
$$\pi_*[\pi^*D] = [D].$$
Donc d'apr\`es ce qui pr\'ec\`ede, 
$$\gamma([D]) = \pi_*\gamma([\pi^*D]).$$
De plus, $\pi^*D$ est aussi principal.

On suppose donc $X$ normale de dimension $n$. On sait alors qu'il existe $T$ un sous-sch\'ema ferm\'e de codimension sup\'erieure ou \'egale \`a deux tel que $U=X-T$ et $Z = [D]\cap U$ soient lisses. On a le diagramme commutatif :
$$\xymatrix{ Z \ar[r]^{i} \ar[d]_{\beta} & U \ar[d]_{\alpha} \\ [D] \ar[r]_{j} & X.}$$
Les axiomes des th\'eories de Poincar\'e et la nullit\'e de la classe d'un diviseur principal donnent
$$\begin{array}{ccc}
\alpha^*\gamma([D]) & = & \alpha^*j_*\eta_{[D]} \\
  & = & i_*\beta^*\eta_{[D]} \\
  & = & i_* \eta_Z \\
  & = & 0.
\end{array}$$

Maintenant en utilisant la suite exacte longue :
$$\cdots H_{2n-2}^{rig}(T) \to  H_{2n-2}^{rig}(X) \surfleche{\alpha^*}   H_{2n-2}^{rig}(U) \to \cdots$$
et sachant que $H_{2n-2}^{rig}(T)=0$, on voit que $\alpha^* : H_{2n-2}^{rig}(X) \to  H_{2n-2}^{rig}(U)$ est injective donc
$$\gamma([D]) = 0.$$

\findem

\begin{prop}
\label{324}
Soient $j : U \hookrightarrow X$ une immersion ouverte et $x \in A_\bullet(X)$. On a :
$$\gamma_Uj^*(x) = j^*\gamma_X(x).$$
\end{prop}

\section{Intersection de cycles}
Nous allons voir que les classes de cycles d\'efinies pr\'ec\'edemment sont compatibles \`a l'intersection. Nous allons dans un premier temps nous limiter au cas o\`u la vari\'et\'e ambiante est lisse avant de traiter le cas g\'en\'eral.

\subsection{Classes de cycles sur les vari\'et\'es lisses}
Soit $X$ une vari\'et\'e lisse de dimension $n$. On se donne deux sous-sch\'emas ferm\'es int\`egres $Y$ et $Z$ de codimension $r$ et $s$ respectivement. Rappelons que l'on  dit que $Y$ et $Z$ s'intersectent {\it proprement} si toutes les composantes de $Y \cap Z$ sont de codimension $r+s$ dans $X$. On peut alors d\'efinir le cycle d'intersection \cite{Fult1} $Y.Z$ qui est un cycle dont le support est $Y \cap Z$.

La vari\'et\'e $X$ \'etant lisse, les classes fondamentales de $Y$ et $Z$ se construisent respectivement dans les groupes $H^{2r}_{Y,rig}(X)$ et $H^{2s}_{Z,rig}(X)$.

\begin{prop}
\label{produit}
Avec les notations pr\'ec\'edentes, si $Y$ et $Z$ s'intersectent proprement on a dans $H^{2(r+s)}_{Y\cap Z, rig}(X)$ :
$$\eta_Y\cup \eta_Z = \eta_{Y.Z}.$$
\end{prop}

\dem Par r\'eduction \`a la diagonale \cite[II.4.2.19]{Gros3}, on peut supposer que $Y$ est lisse. A partir de l\`a, l'immersion $Y \hookrightarrow X$ \'etant r\'eguli\`ere, on se r\'eduit classiquement par une r\'ecurrence au cas o\`u $Y$ est de codimension $1$.

Enon\c{c}ons maintenant un lemme :
\begin{lemme}
Soient $f:X' \to X$ un morphisme de vari\'et\'es int\`egres avec $X$ lisse et $n'$ la dimension de $X'$. Pour tout diviseur de Cartier $D$ sur $X$ tel que $f(X') \nsubseteq |D|$ on a dans $H^{rig}_{2n'-2}(X')$ :
$$\eta_{X'} \cap f^* c_{1,D} = \eta_{[f^*D]}.$$
\end{lemme}

\dem Elle est identique \`a celle du lemme analogue de \cite[7.7.2]{Hart2}.
\findem

On applique le lemme au morphisme $f : Z \hookrightarrow X$, en prenant $Y$ comme diviseur de Cartier. On a alors :
$$\eta_Z \cap f^* c_{1,Y} = \eta_{[Y.Z]}.$$
Il ne reste plus qu'\`a utiliser la formule de projection pour conclure.
\findem

Si $X$ est lisse, le groupe d'homologie de Chow est isomorphe \`a l'anneau de Chow classique. On notera alors aussi 
$$\gamma : A^\bullet X \to H^\bullet_{rig}(X).$$


\begin{prop}
Si $X$ est une vari\'et\'e quasi-projective, le morphisme $\gamma$ est un morphisme d'anneaux.
\end{prop}

\dem Cet \'enonc\'e d\'ecoule de la proposition ci-dessus. En effet si on se donne deux cycles dans $X$, la vari\'et\'e $X$ \'etant quasi-projective, le {\it moving lemma} dit que l'on peut trouver des cycles rationnellements \'equivalents se coupant proprement.
\findem

On sait de plus que pour tout morphisme $f : X' \to X$ entre deux vari\'et\'es lisses quasi-projectives, il existe un morphisme de fonctorialit\'e $f^* : A^\bullet(X) \to A^\bullet(X')$ d\'efini de la mani\`ere suivante. Pour tout cycle $Z$ de $X'$ on choisit un cycle $Z'$ rationnellement \'equivalent \`a $Z$ tel que $p^{-1}Z'$ s'intersecte proprement avec $\Gamma_f$ o\`u $p : X\times X' \to X$ est la projection et $\Gamma_f$ le graphe de $f$. En notant $p': X\times X' \to X'$ l'autre projection, on pose
$$f^*(Z) = p'_*(\Gamma_f.p^{-1}(Z')).$$

\begin{prop}
\label{doli}
Soit $f : X' \to X$ un morphisme entre deux vari\'et\'es lisses quasi-projectives. On a pour tout $z \in A^\bullet(X)$ :
$$f^*\gamma_X(z) = \gamma_{X'}f^*z,$$
o\`u $\gamma_X$ et $\gamma_{X'}$ d\'esignent les applications classes de cycles.
\end{prop}

\dem On suppose d'abord que $f$ est plat. On peut supposer que $z$ est la classe d'un sous-sch\'ema ferm\'e. Par excision on peut de plus supposer que $Z$ est lisse. On se ram\`ene alors par une recurrence au cas de codimension $1$. Ce dernier d\'ecoule des propositions \ref{galois} et \ref{prop331}. 

Pour $f$ quelconque, on peut aussi supposer que $z$ est la classe d'un sous-sch\'ema ferm\'e.
On a alors :
$$\begin{array}{ccc}
f^*\gamma_X(Z) & = & \Gamma_f^*p^*\gamma_X(Z)\\ & = &  \Gamma_f^*\gamma_{X'\times X}(p^*Z) \\
  & = & p'_*(\Gamma_f)_*\Gamma_f^*\gamma_{X'\times X}(p^*Z) \\
& = & p'_* (\gamma_{X'\times X}(\Gamma_f).\gamma_{X'\times X}(p^*Z)) \\
& = & \gamma_{X'}(p'_*(\Gamma_f.p^*Z)) \\
& = & \gamma_{X'}(f^*Z).
\end{array}$$
\findem

\subsection{Cas g\'en\'eral}

On va commencer par rappeler la d\'efinition de l'anneau de cohomologie de Chow pour une vari\'et\'e \'eventuellement singuli\`ere \cite{Fult2}.
Soit $X$ une vari\'et\'e, on note $\CCC(X)$ la cat\'egorie dont les objets sont les morphismes $f : X \to Y$ o\`u $Y$ est une vari\'et\'e lisse. On notera $(Y,f)$ un tel objet. Un morphisme de $(Y,f)$ dans $(Y',f')$ est la donn\'e d'une application $g : Y \to Y'$ telle que $f'=g\circ f.$ On pose alors
$$A^\bullet X = \limind{\CCC(X)} A^\bullet Y.$$

On peut alors construire $\gamma : A^\bullet X \to H^\bullet_{rig}(X).$ En effet soit $(Y,f)$ un objet de $\CCC(X)$, on a 
$$A^\bullet Y \surfleche{\gamma} H^\bullet_{rig}(Y) \surfleche{f^*} H^\bullet_{rig}(X).$$
La proposition \ref{doli} nous montre alors que ces fl\`eches d\'efinissent bien le morphisme voulu. 

\subsection{Cap-produit et classes de cycles}

Il existe un accouplement entre l'homologie et la cohomologie de Chow \cite[3.1]{Fult2} : 
$$A_pX \otimes A^q X \surfleche{\cap} A_{p-q}X.$$
Il est d\'efini de la mani\`ere suivante :
pour tout $x \in A_pX$ et $y \in A^qX$, on choisit $f : X \to Y$ avec $Y$ lisse et $\tilde{y} \in A^qY$ pour repr\'esenter $y$. Alors on pose
$$x\cap y = x \bullet_f \tilde{y},$$  
o\`u $x\bullet_f\tilde{y}$ est l'intersection d\'efinie dans \cite[7]{Ser2}.
On aura pris soin de choisir des cycles s'intersectant proprement.

On a alors :
\begin{prop} 
\label{cap}
Soient $X$ une vari\'et\'e, $x \in A_p X$ et $y \in A^q X$. On a
$$\gamma(x \cap y) = \gamma(x) \cap \gamma(y).$$

\end{prop}

\dem On choisit un rep\'esentant de $y$, \`a savoir un morphisme $f : X \to Y$ o\`u $Y$ est une vari\'et\'e lisse et un cycle $y$ dans $A^q(Y)$.
On doit montrer que
\begin{equation}
\label{Okey}
\gamma_X(x\bullet_f y) = \gamma_X(x) \cap f^*\gamma_Y(y).
\end{equation}
On peut donc se restreindre par lin\'earit\'e au cas o\`u $x$ et $y$ sont repr\'esent\'es par des sous-ch\'emas ferm\'es int\`egres $Z_X$ et $Z_Y$. On se donne de plus, une immersion ferm\'ee de $X$ dans une vari\'et\'e lisse $V$. On a donc le diagramme commutatif suivant :
$$\xymatrix{ X \ar[d]_f \ar[dr]^{\pi}  &   \\ Y & Y\times V. \ar[l]^{p} }$$
En utilisant la proposition \ref{doli}, l'\'egalit\'e (\ref{Okey}) se ram\`ene a :
$$\gamma_X(Z_X\bullet_{\pi} p^*(Z_Y)) = \gamma_X(Z_X) \cap \pi^*\gamma_{Y\times V}(p^*Z_Y).$$
Cette derni\`ere \'egalit\'e est une cons\'equence de la proposition \ref{produit}.

\findem

\section{Applications}
Nous allons maintenant exhiber quelques cons\'equences du formalisme expos\'e pr\'ec\'edemment

\subsection{Th{\'e}or{\`e}me de Riemann-Roch}
Pour commencer, nous allons rappeler la d\'efinition du polyn\^ome de Chern et du polyn\^ome de Todd \cite{SGA6, Hirz}. On se place dans $\QM[[T_i]]_{i\in \NM}$ des s\'eries formelles \`a nombre d\'enombrable de variables. On note $C_i$ les fonctions sym\'etriques des variables. Le polyn\^ome de Chern est d\'efini par :
$${\rm ch}(C_i) = \sum_{i \in \NM} \exp(T_i),$$
o\`u $\exp x = \sum_{i=0}^\infty x^n/n!$.
De m\^eme le polyn\^ome de Todd est d\'efini par :
$${\rm td}(C_i) = \prod_{i\in \NM} \frac{T_i}{1-\exp(-T_i)}.$$

Soient $X$ une $k$-vari\'et\'e et $\ECC$ un faisceau localement libre de rang $r$ sur $X$. On d\'efinit le polyn\^ome de Chern et la classe de Todd de $\ECC$ par
$${\rm ch}(\ECC):= {\rm ch}(c_1(\ECC), \ldots, c_r(\ECC), 0, \ldots ) \hbox{ et }  {\rm td}(\ECC):= {\rm td}(c_1(\ECC), \ldots, c_r(\ECC), 0, \ldots ).$$

Si $X$ est lisse, on appelle classe de Todd de $X$ et on note ${\rm Td}(X)$ la classe de Todd du fibr\'e tangent.

On peut maintenant \'enonc\'e un th\'eor\`eme de Riemman-Roch en cohomologie rigide.

\begin{theo}
Il existe une transformation naturelle $\tau : K_0 \to H_{\bullet}^{rig}$ de foncteurs covariants de la cat\'egorie des $k$-vari\'et\'es quasi-projectives avec des morphismes propres dans la cat\'egorie des groupes ab\'eliens v\'erifiant :
\begin{enumerate}
\item Pour tout $X$, le diagramme 
$$\xymatrix{K^0X \otimes K_0X \ar[r]^{\otimes} \ar[d]_{{\rm ch} \otimes \tau} & K_0X \ar[d]^{\tau} \\ H^{\bullet}_{rig}(X/K) \otimes H_\bullet^{rig}(X/K) \ar[r]^-{\cap} & H_\bullet^{rig}(X/K)}$$
est commutatif.
\item Si $X$ est une $k$-vari\'et\'e lisse,
$$\tau(\OCC_X) = {\rm Td}(X).$$
\item Si $j : U \to X$ est une immersion ouverte de $k$-vari\'et\'es, le diagramme
$$\xymatrix{K_0X \ar[r]^-{\tau} \ar[d]_{j^*} & H_\bullet^{rig}(X/K) \ar[d]^{j^*} \\ K_0 U \ar[r]^-{\tau} & H_\bullet^{rig}(U/K)}$$
est commutatif.
\end{enumerate}

\end{theo}

\dem Le th\'eor\`eme de Riemman-Roch \`a valeur dans l'homologie de Chow \`a coefficient rationnels $A_\bullet(X)_{\QM}$ \cite[II]{BFM}, nous donne l'existence de $\tau_A :   K_0 \to (A_\bullet)_{\QM}$ v\'erifiant les trois hypoth\`eses. On pose alors
$$\tau = \gamma \circ \tau_A.$$
Maintenant, on note $c_{A}$ la th{\'e}orie des classes de Chern {\`a} valeurs dans $A^\bullet(X)_{\QM}$ \cite{SGA6}.

L'unicit{\'e} des classes de Chern {\`a} valeurs dans la
cohomologie rigide nous donne 
$$\gamma \circ c_{univ} = c_{rig}.$$
On a alors pour tout $X$ et tout $\ECC$ faisceau localement libre sur $X$,
$${\rm ch}(\ECC) = \gamma({\rm ch}_A(\ECC)) \hbox{ et } {\rm td}(\ECC) = \gamma({\rm td}_A(\ECC)),$$
o{\`u}  ${\rm ch}_A$ et ${\rm td}_A$ d{\'e}signent respectivement le caract{\`e}re de Chern et la classe de Todd calcul{\'e}s en utilisant $c_{A}$.
Le diagramme commutatif suivant : 

$$\xymatrix{K^0X \otimes K_0X \ar@/^2.5pc/[rrr]^{{\rm ch} \otimes \tau}  \ar[rr]^-{{\rm ch}_A \otimes \tau_A} \ar[d]_{\otimes} &  &
  A^\bullet(X)_{\QM} \otimes A_\bullet(X)_{\QM} \ar[d]^{\cap}  \ar[r]^-{\gamma \otimes \gamma} & H^\bullet_{rig}(X/K)\otimes H_\bullet^{rig}(X/K) \ar[d]^{\cap}\\
 K_0(X) \ar[rr]^{\tau_A}  \ar@/_2.5pc/[rrr]_{\tau}&  &
  A^\bullet(X) \ar[r]^{\gamma} & H_\bullet^{rig}(X/K)}$$
o\`u le carr\'e de droite est commutatif par \ref{cap} nous donne la proposition $1$.

Pour la proposition $2$, on a :
$$\begin{array}{ccc}
\tau(\OCC_X) & = & \gamma({\rm Td}_A(X) \cap [X]) \\
             & = & {\rm Td}(X) \cap \eta_X \\
             & = & {\rm Td}(X).
\end{array}$$
En effet sur $X$ lisse, $\eta_X$ est l'unit\'e de l'anneau de cohomologie.

La proposition $3$ d\'ecoule de la proposition \ref{324}.

\findem

\begin{coro}
Soit $f: X \to Y$ un morphisme propre entre deux $k$-vari\'et\'es lisses quasi-projectives. On a pour tout $\ECC$ dans $K^0X$ :
$$f_*({\rm ch}(\ECC).{\rm Td}(X)) = {\rm ch}(f_*(\ECC)).{\rm Td}(Y).$$
\end{coro}

\subsection{Formule de self-intersection}
Soit $X$ une $k$-vari{\'e}t{\'e} lisse et $Y$ un sous sch{\'e}ma ferm{\'e} lisse de
codimension $d$. On note $\JCC$
l'id{\'e}al de $Y$ dans $X$ et $\NCC=\JCC/\JCC^2$.
\begin{theo}
Pour tout $y \in H^*_{rig}(Y)$ on a :
$$i^*i_*(y) = y.c_d(\check{\NCC}).$$
\end{theo}
\dem
On va commencer par un lemme concernant la cohomologie rigide d'un
{\'e}clat{\'e} sachant que ce th{\`e}me sera repris plus pr{\'e}cis{\'e}ment plus
tard. Fixons les notations. On
consid{\`e}re $f : X' \to X$ l'{\'e}clat{\'e} de $X$ le long de $Y$. On a
alors le diagramme cart{\'e}sien suivant. 
$$\xymatrix{Y' \ar[r]^{j} \ar[d]_{g} & X' \ar[d]^f \\
            Y  \ar[r]_{i}    & X }$$

\begin{lemme}
En notant :
$$f^* : H^*_{Y, rig}(X) \to H^*_{Y', rig}(X')$$
et :
$$f_* : H^*_{Y', rig}(X') \to H^*_{Y, rig}(X)$$
on a l'{\'e}galit{\'e} suivante :
$$f_*\circ f^* = Id.$$
De plus $f_*$ est de degr{\'e} $0$.
\end{lemme}

\dem 
Pour tout $x \in H^*_{Y, rig}(X)$ on a :
$$\begin{array}{ccc}
f_*\circ f^*(x) & = & f_*(1.f^*(x)) \\
             & = & f_*(1).x
\end{array}$$

On est donc ramen{\'e} {\`a} montrer que $f_*(1) = 1$. Il suffit de voir
que c'est vrai sur l'ouvert $X-Y$, ce qui est {\'e}vident car $f$ induit
un isomorphisme de $X'-Y'$ sur $X-Y$. 
\findem

Ce lemme {\'e}tant d{\'e}montr{\'e} la d{\'e}monstration est similiaire {\`a} \cite[p 301]{SGA5}.

\rqe pour les vari\'et\'es quasi-projectives, la formule de self-intersection peut aussi se d\'eduire de la formule dans les groupes de Chow.
\subsection{Cohomologie d'un {\'e}clat{\'e}}
On garde les notations pr{\'e}c{\'e}dentes.
Pour tout $n \in \NM$ on consid{\`e}re l'application :
$$\psi_n : H^n_{rig}(X)\oplus H^{n-2}_{rig}(Y') \to H^n_{rig}(X')
\oplus H^{n-2d}_{rig}(Y)$$
$$(resp.  \psi_{Y,n} : H^n_{Y,rig}(X)\oplus H^{n-2}_{rig}(Y') \to H^n_{Y',rig}(X')
\oplus H^{n-2d}_{rig}(Y))$$
d{\'e}finie par :
$$\psi_n = \left( \begin{array}{cc} f^*  & j_* \\ 0 & g_* \end{array}
\right).$$

\begin{prop}
Avec les notations ci-dessus $\psi_n$ (resp. $\psi_{Y,n}$) est un
isomorphisme.
\end{prop}

\dem
Nous allons d{\'e}montrer que $\psi_{Y,n}$ est un isomorphisme en exhibant
un inverse. Pour commencer on rappelle un lemme de \cite{SGA5} dont la
d{\'e}monstration se recopie dans notre cas.
\begin{lemme}
On note $\FCC$ le $\OCC_{Y'}$-module d{\'e}fini par :
$$0 \to \FCC \to g^*(\NCC) \to \OCC_{Y'}(1) \to 0.$$
Alors pour tout $n\geq 0$ et tout $y \in H^n_{rig}(Y)$ on a :
$$f^*(i_*(y))=j_*(g^*(y)c_{d-1}(\check{\FCC})).$$ 
\end{lemme}

On pose alors :
$$\gamma = j_*^{-1}\circ (f^*\circ f_* - Id).$$
Le lemme pr{\'e}c{\'e}dent nous permet d'affirmer que si on pose 
$$\mu_n = \left( \begin{array}{cc} f_* & -i_* \\ -\gamma &
    c_{d-1}(\FCC)g^* \end{array}\right)$$
alors on a :
$$\mu_n\circ \psi_{Y,n} =Id.$$

On a donc d{\'e}montr{\'e} que pour tout $n$, $\psi_{Y,n}$ est un
isomorphisme.

Passons maintenant {\`a} $\psi_n$.  En remarquant que la
cohomologie de $X-Y$ est la m{\^e}me que celle de $X'-Y'$, les suites exactes
longues d'excisions nous donnent le diagramme commutatif suivant :
$$\xymatrix{ H^{n-1}_{rig}(X-Y) \ar[rr]^{\sim} \ar[d] & &
  H^{n-1}_{rig}(X'-Y') \ar[d] \\ 
H^n_{Y, rig}(X) \oplus H^{n-2}_{Y', rig}(X')
\ar[rr]^{\sim}_{\psi_{Y,n}} \ar[d] &  & H^n_{Y', rig}(X') \oplus
H^{n-2d}_{Y, rig}(X) \ar[d] \\
H^n_{rig}(X) \oplus H^{n-2}_{rig}(X')
\ar[rr]_{\psi_{n}} \ar[d] & & H^n_{rig}(X') \oplus
H^{n-2d}_{rig}(X) \ar[d] \\
H^{n}_{rig}(X-Y) \ar[rr]^{\sim} \ar[d] & &
  H^{n}_{rig}(X'-Y') \ar[d] \\ 
H^{n+1}_{Y, rig}(X) \oplus H^{n+1-2}_{Y', rig}(X')
\ar[rr]^{\sim}_{\psi_{Y,n+1}}&  & H^{n+1}_{Y', rig}(X') \oplus
H^{n+1-2d}_{Y, rig}(X). }$$

Le lemme des cinq nous donne alors la bijectivit{\'e} de $\psi_n$.

\findem

\bibliography{bibthese}

\end{document}